\documentclass[12pt,reqno]{amsart}
\usepackage{amsmath}
\usepackage{amsthm}
\usepackage{graphicx}
\usepackage{hyperref}
\usepackage{mathrsfs}
\usepackage{cite}
\usepackage{dsfont}
\usepackage{color}
\usepackage{bm}
\usepackage{amsfonts}
\usepackage{amssymb}
\usepackage{enumerate}
\usepackage{array}
\usepackage{cases}
\usepackage{soul}
\usepackage[margin=3cm]{geometry}
\usepackage[normalem]{ulem}

\makeatletter
\def\@settitle{\begin{center}%
		\baselineskip14\p@\relax
		\bfseries
		\uppercasenonmath\@title
		\@title
		\ifx\@subtitle\@empty\else
		\\[1ex]\uppercasenonmath\@subtitle
		\footnotesize\mdseries\@subtitle
		\fi
	\end{center}%
}
\def\subtitle#1{\gdef\@subtitle{#1}}
\def\@subtitle{}
\makeatother

\numberwithin{equation}{section}

\newtheorem{theorem}{Theorem}[section]
\newtheorem{lemma}[theorem]{Lemma}
\newtheorem{proposition}[theorem]{Proposition}

\newtheorem{definition}[theorem]{Definition}

\definecolor{darkgreen}{rgb}{0,0.5,0}
\definecolor{darkblue}{rgb}{0,0,0.7}
\definecolor{darkred}{rgb}{0.9,0.1,0.1}
\definecolor{lightblue}{rgb}{0,0.51,1}

\newcommand{\R}{{\mathbb R}}

\newcommand{\T}{{\mathcal T}}

\newcommand{\pa}{\partial}
\newcommand{\dt}{{\partial_t}}

\newcommand{\dxi}{\partial_{x_i}}
\newcommand{\dxj}{\partial_{x_j}}

\DeclareMathOperator{\divv}{div}

\def\R{{\mathbb R}}

\definecolor{darkspringgreen}{rgb}{0.09, 0.45, 0.27}

\title[]{Regularity estimates for diffusion semigroups on weighted Sobolev spaces}
\author[]{M. Hauray, Y. V. Vuong}

\begin{document}
	
	\maketitle

	\begin{abstract}
		In this work, we consider a class of second order uniformly elliptic operators with smooth and bounded coefficients. We provide some estimates on the norm of the semigroup generated by these operators acting on weighted Sobolev spaces, where the weight satisfies some specific conditions. Our proof relies on a classical bound for the derivatives of fundamental solution to parabolic equations.

	\end{abstract}

	\section{Introduction}
	
	We consider a second order differential operator $A$ on $\R^d$ given in divergence form by
	\[A\phi=- \sum_{i,j=1}^{d} \dxi\big(a_{ij}(x)\dxj \phi\big),\] where the coefficients $a_{ij}$ are symmetric, bounded and smooth enough (to be precised later) and satisfy an uniform ellipticity condition, i.e.
	\[\sum_{i,j=1}^{d} a_{ij}(x)\xi_i\xi_j\ge \lambda |\xi|^2, \quad \forall x, \xi \in \R^d, \]
	for some positive constant $\lambda$. Under these assumptions, $A$ is a self-adjoint and positive operator on $L^2(\R^d)$. Thanks to Hille-Yoshida theorem, $A$ generates the strongly continuous semigroup $\T_t=e^{-tA}$ on $L^2(\R^d)$.
	
	\smallskip
	Denoting
	\begin{equation*}
		u(t,x)=	\T_t \phi (x), \qquad t\ge 0,\; x\in \R^d,
	\end{equation*}
	we obtain the unique solution to the Cauchy problem
	
	\begin{equation}\label{eq:parabolic_pde}
		\begin{cases}
			\dt u(t,x)&=-A u(t,x), \qquad t>0, \; x\in \R^d,\\
			u(0,x)&=\phi (x), \qquad \qquad \qquad \quad \, x\in \R^d.
		\end{cases}
	\end{equation}
	%admits a unique classical solution with $\phi\in L^2(\R^d)$, i.e., there exists a unique bounded and continuous function $u: [0, +\infty) \times \R^d \to \R$ which is once continuously diﬀerentiable with respect to time and twice continuously diﬀerentiable with respect to the space variables in $(0, +\infty) \times \R^d$, and such that $u(0, x)=\phi(x)$. This allows us to define a semigroup $(\T_t)_{t\ge 0}$ of linear operators in $C(\R^d)$

	\medskip
	
	Our main result in this article concerns the regularity of the semigroup $\T_t$ in a class of weighted Sobolev spaces. To state properly our main theorem, we need to introduce some notations related to the weights and the weighted spaces.
	
	\medskip
	Let $w:\R^d \to \R$ be a nonnegative $C^1$ function that satisfies
	\begin{align}
		|\nabla w(x)|&\le C_1 w(x), \qquad \forall x\in \R^d, \label{eq:derivative_of_w}\\
		w(x)w(y)^{-1 }&\le C_2\exp\big(C_3|x-y|\big), \qquad \forall x, y\in \R^d, \label{eq:w(x)w(y)-condition}
	\end{align} for some positive constants $C_1, C_2, C_3$.
	
	\medskip
	We will give two important examples of weights satisfying the above assumptions.
	
	\begin{proposition}\label{prop:weight_examples}
		The two families of weights (polynomial and exponential) satisfy conditions \eqref{eq:derivative_of_w} and \eqref{eq:w(x)w(y)-condition}:
		\begin{itemize}
			\item[$i)$]  $\displaystyle w_{\alpha}(x)=(1+|x|^2)^{\alpha}$  with $\alpha \in\R$,
			\item[$ii)$]  $\displaystyle w_{\lambda}(x)=\exp(\lambda\sqrt{1+|x|^2})$ with $\lambda \in\R$.
		\end{itemize}
	\end{proposition}
	
	The proof of Proposition~\ref{prop:weight_examples} is given in Appendix \ref{sec:example}.
	
	\begin{definition}
		The weighted Sobolev space $H^k_{w}$ , consists of all locally summable, $k$ times weakly differentiable  function $\phi:\R^d\to \R$ for which the norm $$\|\phi\|_{H^k_{w}}:=\left(\sum_{|m|\le k}\int |\nabla^m \phi(x)|^2 w(x)dx\right)^{1/2}$$
		is finite.
	\end{definition}
	
	The motivation for studying this kind of weighted Sobolev spaces can be found in \cite{Hauray-Pardoux-Vuong22}, where the authors use the polynomial weight $w_{\alpha}(x)=(1+|x|^2)^{-\alpha},\; \alpha>0$ and the results obtained in this paper to prove a central limit theorem for an mean field interacting system in epidemiology. The estimates \eqref{eq:weight:semigroup_bound_first_estimate} and \eqref{eq:weight:semigroup_bound_second_estimate} are crucial in the proof of central limit theorem in that paper.
	
	\medskip
	The main result of this paper is the following.
	
	\begin{theorem}\label{prop:weight:regularity_of_semigroups}
		Let $k\ge 0$ and assume that $a_{ij}\in C^{2k+1}_b(\R^d)$. Let $(\T_t)_{t\ge 0}$ be the semigroup generated by $A$. For any $T\ge 0$, there exists a constant $C_T>0$ depends only on $T, d, k, \|a\|_{H^{2k+1}_{w}}$ such that for any $t\in [0,T]$, the following holds true
		\begin{enumerate}
			\item
			\begin{equation}\label{eq:weight:semigroup_bound_first_estimate}
				\|\T_t\phi\|_{H^{k}_{w}}\le C_T\|\phi\|_{H^{k}_{w}}.
			\end{equation}
			\item
			\begin{equation}\label{eq:weight:semigroup_bound_second_estimate}
				\|\nabla_x\T_t\phi\|_{H^{k}_{w}}\le C_T\left(1+\frac{1}{\sqrt{t}}\right)\|\phi\|_{H^{k}_{w}}.
			\end{equation}
		\end{enumerate}
	\end{theorem}
	
	%where $H^k_{\alpha}$ are weighted Sobolev spaces which are endowed with the $\|\cdot\|_{H^k_{\alpha}}$ norms, defined by $$\|\phi\|_{H^k_{\alpha}}:=\left(\sum_{|m|\le k}\int |\nabla^m \phi(x)|^2 m_{\alpha}(x)dx\right)^{1/2},$$
	%for any $\phi\in H^k_{\alpha}$, and $$m_{\alpha}(x)=\frac{1}{(1+|x|^2)^{\alpha}}.$$
	%
	%An important property of this weight is,
	%\begin{align}\label{eq:derivative_of_w}
	%	|\nabla m_{\alpha}(x)|=\frac{2\alpha|x|}{(1+|x|^2)^{\alpha+1}}\le 2\alpha\, m_{\alpha}(x).
	%\end{align}
	%In fact, for any weight $m$ which has the property $|\nabla m(x)|\le C m(x)$ for some positive constant $C$, all the results obtained in this article still hold true.
	%
	%An example for the application as well as the motivation for the study of this kind of estimate can be found in \cite{Hauray-Pardoux-Vuong22}, where the authors use the estimates \eqref{eq:weight:semigroup_bound_first_estimate}, \eqref{eq:weight:semigroup_bound_second_estimate} obtained in the present article with the choice of weight $m_{\alpha}(x)=(1+|x|^2)^{-\alpha}$ for an evolution equation in epidemiology.
	
	\medskip
	
	Recently, the study of uniformly elliptic operators with unbounded coefficients has been an active field of research in PDEs theory. The motivations for this study derive from the fact that operators with unbounded coefficients naturally emerge in the theory of Markov diffusion processes. We notice that in this case, even the differentiability of the semigroup is nontrivial since singular diffusion coefficients (in case of degenerate noises) can lead to loss of regularity, see
	Theorem~1.2 in \cite{Hairer15}. Under suitable assumptions on the growth of coefficients, the authors in \cite{Lunardi98} provide a pointwise estimate for the derivatives of semigroup up to order three, i.e. there exists $\delta\in \R$ such that 
	\begin{equation}
		\|\T_t\|_{L(C(\R^d),C^k(\R^d))}\le \frac{Ce^{\delta t}}{t^{k/2}}, \quad \text{for } k=1,2,3.
	\end{equation}
	See also \cite{Bertoldi-Lorenzi05}, \cite{Lorenzi05} for the other versions with some relaxations on the assumptions.
	
	\medskip
	There are two main approach to study the problem of estimating the derivatives of $\T_t\phi$: the analytic methods \cite{Lunardi97}, \cite{Lunardi98}, \cite{Prato99} and probabilistic methods \cite{Cerrai98}, \cite{Priola-Wang06}, \cite{Wang97}. In the analytical approach, they consider the analytic extension of the $C_0$-semigroup ("analytic" means that $(\T_t)_{t\ge 0}$ is analytic w.r.t. $t$ in some sector containing the non-negative real axis, see Definition $5.1$ in \cite{Pazy83}), and study the connection between the analytic semigroup and its generator. A fundamental and important property of the analytic semigroups is the following: If $(\T_t)_{t\ge 0}$ is an analytic semigroup generated by $A$, then for any $t>0$, the operator $A^{p}\T_t$ ($p>0$) is bounded in $L^2(\R^d)$ and
	\begin{equation}
		\|A^{p}\T_t\| \le C_{p} t^{-p},
	\end{equation}
	for some constant $C_{p} >0$, and $\|\cdot\|$ here is the operator norm on $L^2(\R^d)$ (see Theorem $6.13$ in \cite{Pazy83}). In contrast, the probabilistic approach relies on coupling techniques and the Bismut–Elworthy-Li formula, which gives an explicit expression for the derivatives of $\T_t\phi$ in term of the derivatives of $\phi$ and the solution to the associated SDE (see \cite{Bismut84} and the extended version of this formula \cite{Elworthy-Li94}).
	
	%In \cite{Hutzenthaler-Pieper21}, from a probabilistic viewpoint, the authors provide some regularity estimates for a specific type of Markov transition semigroup with constant independent of the dimension, and seems to be a new . It allows to analyze the high dimensional equations.
	
	%the solution to the backward Kolmogorov equation can be represented as an
	%expectation of a functional of the solution of an SDE (see \cite{Bismut84} and the extended version of this formula \cite{Elworthy-Li94}).
	
	\medskip
	Although many works have been done to study this kind of estimate on spatial derivatives of diffusion semigroups in both bounded and unbounded cases, to the best of our knowledge, there is no result dealing with weighted norms. In fact, it also seems difficult to extend the previous methods in a weighted case.
	
	\medskip
	Our approach is more classical and we only deals with bounded and smooth coefficients. First, we take advantage of the results on the differentiability of the fundamental solution of parabolic equations to prove Theorem~\ref{prop:weight:regularity_of_semigroups} in the case of classical Sobolev spaces without weights. After that, we extend the results obtained in Section~\ref{sec:case_without_weights} to the weighted case.

	\section{The case without weights}\label{sec:case_without_weights}
	
	Along the study of regularity estimates of the semigroup $(\T_t)_{t\ge 0}$ in this section, we consider the classical Sobolev spaces $H^k$ ($H^k_w$ with $w=1$) which are endowed with the $\|\cdot\|_{H^k}$ norms, $\|\phi\|_{H^k}:=\left(\sum_{|m|\le k}\int |\partial^m \phi(x)|^2 dx\right)^{1/2}$. This case is covered by previous results that could be found in the literature, but we will provide a proof that could be extended to the weighted case.
	
	\medskip
	To prove Theorem \ref{prop:weight:regularity_of_semigroups} in the case without weights, the main ingredient that we use is a classical bound on the fundamental solution to the parabolic equation \eqref{eq:parabolic_pde} (see \cite{Friedman83}, Theorem $2$, Chapter $9$). First, we provide the estimates for the low regularity spaces $H^k$ with $k=0,1$, then by induction, we extend it to higher regularity spaces using the commutation of $A$ and the associated semigroup $\T_t$.
	
	\medskip
	We start with a lemma which roughly states the equivalence between $\|A^k\phi\|_{L^2}$ and $\|\phi\|_{H^{2k}}$ for any $k\ge 1$.
	\begin{lemma}\label{lem:Hk_estimate_for_induction} 
		Let $k\ge 1$ and assume that $a_{ij}\in C^{2k-1}_b(\R^d)$.
		
		\begin{enumerate} 
			\item There exists $C_1,C_2>0$ such that
			\begin{equation}\label{eq:Hk_estimate_for_induction}
				C_1\|\phi\|_{H^{2k}}\le \|A^k\phi\|_{L^2}+\|\phi\|_{H^{2k-1}}\le C_2\|\phi\|_{H^{2k}}.
			\end{equation}
			
			\item If we assume that $a_{ij}\in C^{2k}_b(\R^d)$, then there exists a positive constant $C$ such that
			\begin{equation}\label{eq:lem:Akphi_H1_estimate}
				\|A^k\phi\|_{H^1}\le C\|\phi\|_{H^{2k+1}}.
			\end{equation}
		\end{enumerate}
	\end{lemma}
	
	\begin{proof} $1).$ \textbf{Step 1.} We first prove the inequality on the r.h.s. of \eqref{lem:Hk_estimate_for_induction}. By the definition of the operator $A$, the operator $A^k$ can be represented in the following form
		\begin{align}\label{eq:expansion_of_Akphi}
			A^k\phi =&\sum_{m=1}^{2k}\sum_{l_1,\dots,l_m=1}^d b^m_{l_1\dots l_m}\pa^{m}_{l_1\dots l_m}\phi,
		\end{align}
		where $b^m_{l_1\dots l_m}$ is a sum of the products consisting of $a_{ij}$ and its derivatives up to order $2k-m$, and hence $b^m_{l_1\dots l_m}\in C^{m-1}_b(\R^d)$ thanks to the assumption on $a_{ij}$. Using the above expansion of $A^k$ we have
		\begin{equation}\label{eq:upperbound_A^k_norm_ineq1}
			\begin{aligned}
				\|A^k\phi\|_{L^2}^2=&\big<A^k\phi, A^k\phi\big>\\
				=& \sum_{m, m'=1}^{2k}\sum_{\substack{l_1,\dots,l_m=1\\ l_1',\dots, l_m'=1}}^d \int b^m_{l_1\dots l_m}b^{m'}_{l_1'\dots l_m'}\pa^{m}_{l_1\dots l_m}\phi\,\pa^{m'}_{l_1'\dots l_m'}\phi.
			\end{aligned}
		\end{equation}
		By Cauchy-Schwartz and the fact that all coefficients $b^m_{l_1\dots l_m}$ are bounded,
		\begin{equation}\label{eq:upperbound_A^k_norm_ineq2}
			\int b^m_{l_1\dots l_m}b^{m'}_{l_1'\dots l_m'}\pa^{m}_{l_1\dots l_m}\phi\pa^{m'}_{l_1'\dots l_m'}\phi \le C\|\pa^{m}\phi\|_{L^2}\|\pa^{m'}\phi\|_{L^2}\le C \big(\|\pa^{m}\phi\|_{L^2}^2+\|\pa^{m'}\phi\|_{L^2}^2\big).
		\end{equation}
		Hence \eqref{eq:upperbound_A^k_norm_ineq1} and \eqref{eq:upperbound_A^k_norm_ineq2} implies
		\begin{equation}\label{eq:upperbound_A^k_norm_ineq3}
			\|A^k\phi\|_{L^2}^2\le C\sum_{m=1}^{2k} \|\pa^m \phi\|_{L^2}^2\le C\|\phi\|_{H^{2k}}^2,
		\end{equation}
		and the inequality on the r.h.s. of \eqref{eq:Hk_estimate_for_induction} follows.	
		
		\medskip
		\textbf{Step 2.}
		To prove the inequality on the l.h.s. of \eqref{eq:Hk_estimate_for_induction}, we use induction for smooth function $\phi$. First, we treat the case $k=1$. Using the uniform ellipticity condition we have
		\begin{equation}\label{eq:lowebound_A^k_techstep1}
			\begin{aligned}
				\left<\phi,A\phi\right>=&-\sum_{i,j=1}^{d}\int \phi\dxi(a_{ij}\dxj\phi)\\
				=&\sum_{i,j=1}^{d}\int (\dxi\phi) a_{ij}\dxj\phi\\
				\ge& \lambda\|\nabla\phi\|_{L^2}^2.
			\end{aligned}
		\end{equation}
		Applying the above estimate to $\nabla \phi$, we get
		\begin{equation}\label{eq:proof:norms_equiv:2nabla_bound}
			\begin{aligned}
				\|\nabla\nabla\phi\|_{L^2}^2 \le& \frac{1}{\lambda}\big<\nabla \phi,A(\nabla\phi)\big>\\
				\le& \frac{1}{\lambda}\Big(\big<\nabla \phi,\nabla(A\phi)\big>-\big<\nabla \phi,[\nabla, A]\phi\big>\Big),
			\end{aligned}
		\end{equation}
		where $[\nabla, A]$ is the commutator between $\nabla$ and $A$, defined by \begin{equation}\label{commutator}
			[\nabla, A]\phi:=\nabla(A\phi)-A(\nabla\phi)=\sum_{i,j=1}^{d} \dxi\big(\nabla a_{ij}\dxj \phi\big).
		\end{equation}
		Noticing that $[\nabla, A]$ is also a second order operator in divergence form and we can get the same estimate as in \eqref{eq:upperbound_A^k_norm_ineq3} with $k=1$,
		\begin{equation}\label{bound_for_commutator}
			\|[\nabla, A]\phi\|_{L^2}\le C \left(\|\nabla ^2 \phi\|_{L^2}+\|\nabla\phi\|_{L^2}\right).
		\end{equation}
		Now by using Cauchy-Schwartz, we deduce from \eqref{eq:proof:norms_equiv:2nabla_bound} that
		\begin{equation}
			\begin{aligned}
				\|\nabla^2\phi\|_{L^2}^2
				\le& \frac{1}{\lambda}\Big(-\big<\Delta \phi,A\phi\big>-\big<\nabla \phi,[\nabla, A]\phi\big>\Big)\\
				\le& C\Big(\|\Delta\phi\|_{L^2}\|A\phi\|_{L^2}+\|\nabla \phi\|_{L^2}\|[\nabla, A]\phi\|_{L^2}\Big)\\
				%		\le& C\Big(\|\nabla ^2\phi\|_{L^2}\|A\phi\|_{L^2}+\|\nabla \phi\|_{L^2}\left(\|\nabla ^2 \phi\|_{L^2}+\|\nabla\phi\|_{L^2}\right)\Big)\\
				\le& C\Big(\|\nabla ^2\phi\|_{L^2}\|A\phi\|_{L^2}+\|\nabla \phi\|_{L^2}\|\nabla ^2 \phi\|_{L^2}+\|\nabla\phi\|_{L^2}^2\Big).
			\end{aligned}
		\end{equation}
		Since $\|\nabla\phi\|_{L^2}^2=\big<\nabla\phi, \nabla\phi\big>=-\big<\Delta\phi, \phi\big>\le \|\Delta \phi\|_{L^2}\|\phi\|_{L^2}\le \|\nabla^2 \phi\|_{L^2}\|\phi\|_{L^2}$, we finally get
		\begin{equation}\label{eq:lowebound_A^k_techstep2}
			\begin{aligned}
				\|\nabla^2\phi\|_{L^2} \le&
				C\Big(\|A\phi\|_{L^2}+\|\nabla \phi\|_{L^2}+\|\phi\|_{L^2}\Big)\\
				\le& C\Big(\|A\phi\|_{L^2}+\|\phi\|_{H^1}\Big).
			\end{aligned}
		\end{equation}
		By iterating the above inequality, we can obtain the estimates for differential operators $A^k$ with higher order $k$. Indeed, assuming that the inequality on the l.h.s. of \eqref{eq:Hk_estimate_for_induction} is true up to $k$, we will prove it holds true for $k+1$. We have
		\begin{equation}\label{eq:lowebound_A^k_ineq1}
			\begin{aligned}
				\|\nabla^{2k+2}\phi\|_{L^2}=&\|\nabla^2(\nabla^{2k}\phi)\|_{L^2}\\
				\le& C\Big(\|A(\nabla^{2k}\phi)\|_{L^2}+\|\nabla^{2k}\phi\|_{H^1}\Big)\\
				\le& C\Big(\|\nabla^{2k}(A\phi)\|_{L^2}+\|[\nabla^{2k},A]\phi\|_{L^2}+\|\nabla^{2k}\phi\|_{H^1}\Big).
			\end{aligned}
		\end{equation}
		Since $[\nabla,A]$ is a second order differential operator, the following expansion
		
		\[
		[\nabla^{2k},A]=\sum_{m=0}^{2k-1} \nabla^m [\nabla,A]\nabla^{2k-1-m}
		\]
		implies that $[\nabla^{2k},A]$ is a differential operator of order $2k+1$. Hence, if $a_{ij}$ belongs to $C_b^{2k+1}(\R^d)$ then $[\nabla^{2k},A]$ is well-defined and
		\begin{equation}\label{eq:bound:commutator_nabla2k_A}
			\|[\nabla^{2k},A]\phi\|_{L^2}\le \|\phi\|_{H^{2k+1}}.
		\end{equation}
		So from \eqref{eq:lowebound_A^k_ineq1} and \eqref{eq:bound:commutator_nabla2k_A} we get
		\begin{equation*}
			\begin{aligned}
				\|\nabla^{2k+2}\phi\|_{L^2}
				\le& C\Big(\|\nabla^{2k}(A\phi)\|_{L^2}+\|\phi\|_{H^{2k+1}}\Big).
			\end{aligned}
		\end{equation*}
		Now by the induction hypothesis, we can conclude that
		\begin{equation}\label{eq:lowebound_A^k_ineq2}
			\begin{aligned}
				\|\nabla^{2k+2}\phi\|_{L^2}
				\le& C\Big(\|A^{k+1}\phi\|_{L^2}+\|A\phi\|_{H^{2k-1}}+\|\phi\|_{H^{2k+1}}\Big)\\
				\le& C\Big(\|A^{k+1}\phi\|_{L^2}+\|\phi\|_{H^{2k+1}}\Big),
			\end{aligned}
		\end{equation}
		which completes the proof of \eqref{eq:Hk_estimate_for_induction}.

		\medskip
		$2).$ The proof of \eqref{eq:lem:Akphi_H1_estimate} is similar to the proof of the inequality on the r.h.s. of \eqref{lem:Hk_estimate_for_induction}. Indeed, from the expansion of $A^k\phi$ \eqref{eq:expansion_of_Akphi} we also have the following
		\begin{align*}
			\pa_{x_l}A^k\phi =&\sum_{m=1}^{2k+1}\sum_{l_1,\dots,l_m=1}^d b^m_{l_1\dots l_m}\pa^{m}_{l_1\dots l_m}\phi, \quad l=1,\ldots, d.
		\end{align*}
		Therefore, using the same estimates \eqref{eq:upperbound_A^k_norm_ineq1}, \eqref{eq:upperbound_A^k_norm_ineq2}, we deduce that
		\begin{equation}\label{eq:Akphi_H1norm_ineq1}
			\|\pa_{x_l}A^k\phi\|_{L^2}^2\le C\sum_{m=1}^{2k+1} \|\pa^m \phi\|_{L^2}^2\le C\|\phi\|_{H^{2k+1}}^2.
		\end{equation}
		Finally, we obtain \eqref{eq:lem:Akphi_H1_estimate} by combining \eqref{eq:upperbound_A^k_norm_ineq3} and \eqref{eq:Akphi_H1norm_ineq1}.
		
	\end{proof}
	
	\medskip
	
	Next, we will use an important result on the differentiability of fundamental solution to parabolic equations to prove regularity of the semigroup generated by $A$ on $L^2$. The following lemma is a consequence of (\cite{Friedman83}, Theorem $2$, Chapter $9$) applied to the parabolic equation \eqref{eq:parabolic_pde}.
	\begin{lemma}\label{lem:differentiability_of_fundamental_solution}
		Let $T>0$. Assume that the matrix $(a)_{ij}$ is uniformly elliptic and all the coefficients $a_{ij}$ belong to $C^1_b(\R^d)$. Then there exists a fundamental solution $K$ of equation \eqref{eq:parabolic_pde} satisfying
		\begin{equation}
			\forall\, t\le T,\; \forall\, x,y\in \R^d, \quad \big|\nabla_x K(t,x,y)\big|\le \frac{C}{t^{(d+1)/2}}\exp\left(-c\frac{|x-y|^2}{t}\right),
		\end{equation}
		where $C, c$ are positive constants.
	\end{lemma}
	
	In fact, Lemma \ref{lem:differentiability_of_fundamental_solution} is a very special case of Theorem $2$ in \cite{Friedman83} with $p=2$, $m=1$, which are respectively the order of the equation and the number of spatial derivatives. In addition, we made a smoothness assumption on the coefficients which is stronger than the one in \cite{Friedman83} (Holder continuity only).
	
	\medskip
	Now, we will use Lemma \ref{lem:Hk_estimate_for_induction} and Lemma \ref{lem:differentiability_of_fundamental_solution} to prove a version of Theorem \ref{prop:weight:regularity_of_semigroups} in the classical Sobolev norms. 
	
	\begin{proposition}\label{prop:regularity_of_semigroups}
		Let $k\ge 0$ and assume that $a_{ij}\in C^{2k+1}_b(\R^d)$. Let $(\T_t)_{t\ge 0}$ be the semigroup generated by $A$. Then, there exists a constant $C>0$ depending only on $d, k, \|a\|_{H^{2k+1}}$ such that for any $t\ge 0$,
		\begin{enumerate}
			\item
			\begin{equation}\label{eq:semigroup_bound_first_estimate}
				\|\T_t\phi\|_{H^{k}}\le C\|\phi\|_{H^{k}}.
			\end{equation}
			\item
			\begin{equation}\label{eq:semigroup_bound_second_estimate}
				\|\nabla\T_t\phi\|_{H^{k}}\le C\left(1+\frac{1}{\sqrt{t}}\right)\|\phi\|_{H^{k}}.
			\end{equation}
		\end{enumerate}
	\end{proposition}
	
	\begin{proof}
		We will prove both \eqref{eq:semigroup_bound_first_estimate} and \eqref{eq:semigroup_bound_second_estimate} simultaneously by induction. At the first step, we will prove \eqref{eq:semigroup_bound_first_estimate} for the cases $k=0$ and $k=1$. The second step is devoted to the cases $k=0,1$ of \eqref{eq:semigroup_bound_second_estimate}. Next, we will use Lemma \eqref{eq:Hk_estimate_for_induction} to perform induction, but the recurrence step is not straightforward. It depends on the parity of the step, i.e. once inequality \eqref{eq:semigroup_bound_first_estimate} is true for $2k$, we prove that it  holds true for $2k+1$, and once it is true for $2k+1$, we prove that it is also true for $2k+2$. In the last step we perform the induction for \eqref{eq:semigroup_bound_second_estimate}.
		
		\medskip
		\paragraph{\textbf{Step 1.} } First, we prove \eqref{eq:semigroup_bound_first_estimate} for $k=0$ and $k=1$.
		\begin{itemize}
			\item \underline{Case $k=0$:}
		\end{itemize}
		
		Since $\T_t\phi$ is a solution of equation \eqref{eq:parabolic_pde}, one has
		
		\begin{equation}\label{eq:Solution_to_parabolic_eq_1}
			\dt\big(\T_t\phi\big)=\divv\big(a\nabla (\T_t\phi)\big).
		\end{equation}
		We multiply both sides of \eqref{eq:Solution_to_parabolic_eq_1} by $\T_t\phi$ and integrate in position
		\begin{align}\label{eq:proof_regularity_semigroup_step1_L^2_integration}
			\frac12\|\T_t\phi\|_{L^2}^2 =& \frac12\|\T_0\phi\|_{L^2}^2+\int_0^t \int_{\R^d}\T_s\phi \divv\big(a\nabla (\T_s\phi)\big)ds \nonumber\\
			=& \frac12\|\phi\|_{L^2}^2-\int_0^t \int_{\R^d}\nabla(\T_s\phi) a\nabla (\T_s\phi)ds.
		\end{align}
		Since the matrix $a$ satisfies the uniform ellipticity condition, we have
		\begin{equation}\label{eq:proof_regularity_semigroup_step1_elliptic_cond}
			\nabla(\T_s\phi) a\nabla (\T_s\phi)\ge \lambda|\nabla (\T_s\phi)|^2.
		\end{equation}
		So \eqref{eq:proof_regularity_semigroup_step1_L^2_integration} and \eqref{eq:proof_regularity_semigroup_step1_elliptic_cond} lead to
		\begin{align}\label{eq:proof_regularity_semigroup_step1_k=0}
			\begin{split}
				\|\T_t\phi\|_{L^2}^2 
				\le& \|\phi\|_{L^2}^2 -2\lambda\int_0^t \|\nabla (\T_s\phi)\|_{L^2}^2ds\\
				\le& \|\phi\|_{L^2}^2,
			\end{split}
		\end{align}
		and \eqref{eq:semigroup_bound_first_estimate} is true for $k=0$.
		
		\medskip
		\begin{itemize}
			\item \underline{Case $k=1$:}
		\end{itemize}
		
		Now for each $l\in \{1,\ldots,d\}$, we take the derivative with respect to $x_l$ on both sides of \eqref{eq:Solution_to_parabolic_eq_1} and obtain the following system
		
		\begin{equation}\label{eq:Solution_to_parabolic_system}
			\dt\psi_l(t) =\divv\big(a\nabla \psi_l(t) \big)+\divv\big((\pa_{x_l} a) \psi(t)\big),\quad l=1,\ldots,d
		\end{equation}
		where $\psi_l(t):=\pa_{x_l}\T_t\phi$ and $\psi=\big(\psi_1,\ldots,\psi_d\big)=\nabla\T_t\phi$.
		
		\smallskip
		Similarly, we can also obtain an appropriate bound for $|\psi(t)|^2$ from equation \eqref{eq:Solution_to_parabolic_system} by multiplying both sides of \eqref{eq:Solution_to_parabolic_system} with $\psi_l$, then integrating in position
		\begin{align}\label{eq:proof_regularity_semigroup_step2_L^2_integration_psi}
			\frac12\|\psi_l(t)\|_{L^2}^2 =& \frac12\|\psi_l(0)\|_{L^2}^2+\int_0^t \int_{\R^d} \psi_l(s) \divv\big(a\nabla \psi_l(s) \big)ds +\int_0^t \int_{\R^d} \psi_l(s) \divv\big((\pa_{x_l} a) \psi(s)\big)ds \nonumber\\
			=& \frac12\|\pa_{x_l}\phi\|_{L^2}^2-\int_0^t \int_{\R^d}\nabla\psi_l(s) a\nabla \psi_l(s)ds -\int_0^t \int_{\R^d}\nabla\psi_l(s) (\pa_{x_l} a) \psi(s)ds.
		\end{align}
		Again we use the uniform ellipticity condition of the matrix $a$ to deduce that 
		\begin{equation}\label{eq:proof_regularity_semigroup_step2_ineq3}
			\nabla\psi_l(s) a\nabla \psi_l(s)\ge  \lambda|\nabla\psi_l(s)|^2 .
		\end{equation}
		On the other hand, all the coefficients $(\pa_{x_l} a_{ij})_{i,j,l}$ are uniformly bounded. So using Cauchy-Schwartz we can quantify the last term by
		\begin{align}\label{eq:proof_regularity_semigroup_step2_ineq4}
			- \int_{\R^d}\nabla\psi_l(s) (\pa_{x_l} a) \psi(s)
			\le& C \int_{\R^d}|\nabla\psi_l(s)| |\psi(s)| \nonumber\\
			\le& \lambda \|\nabla\psi_l(s)\|_{L^2}^2 +\frac{C^2}{4\lambda} \|\psi(s)\|_{L^2}^2 .
		\end{align}
		Combining \eqref{eq:proof_regularity_semigroup_step2_L^2_integration_psi}, \eqref{eq:proof_regularity_semigroup_step2_ineq3} and \eqref{eq:proof_regularity_semigroup_step2_ineq4} we finally obtain
		\begin{align}\label{eq:proof_regularity_semigroup_step2_L^2_estimate_psi}
			\|\psi_l(t)\|_{L^2}^2
			\le& \|\pa_{x_l}\phi\|_{L^2}^2+C\int_0^t \|\psi(s)\|_{L^2}^2 ds,\quad l=1,\ldots,d
		\end{align}
		Now summing up on $l=1,\ldots,d$ and with the help of Gronwall's lemma we deduce that
		
		\begin{equation}\label{eq:proof_regularity_semigroup_step2_L^2_estimate_psi_Gronwall2}
			\|\nabla\T_t\phi\|_{L^2}^2=\|\psi(t)\|_{L^2}^2\le e^{Ct}\|\nabla\phi\|_{L^2}^2,
		\end{equation}
		which completes the proof of \eqref{eq:semigroup_bound_first_estimate} for the case $k=1$.
		
		\medskip
		
		\paragraph{\textbf{Step 2.} }
		Next, we will verify that the inequality \ref{eq:semigroup_bound_second_estimate} is also true for $k=0$ and $k=1$. In other words, we need to prove that
		\begin{equation}\label{eq:proof_regularity_semigroup_step2_ineq1}
			\|\nabla\T_t\phi\|_{L^2}\le C\left(1+\frac{1}{\sqrt{t}}\right)\|\phi\|_{L^2},
		\end{equation}
		and
		\begin{equation}\label{eq:proof_regularity_semigroup_step2_ineq2}
			\|\nabla^2\T_t\phi\|_{L^2}\le C\left(1+\frac{1}{\sqrt{t}}\right)\|\phi\|_{H^1}.
		\end{equation}
		\begin{itemize}
			\item \underline{Case $k=0$:}
		\end{itemize}
		
		Let $K$ be the fundamental solution of equation \eqref{eq:parabolic_pde}. By Lemma \ref{lem:differentiability_of_fundamental_solution} we have
		
		\begin{equation}\label{eq:1_diff_bound}
			\forall\, t\le T, \quad \big|\nabla_x K(t,x,y)\big|\le \frac{C}{t^{(d+1)/2}}\exp\left(-c\frac{|x-y|^2}{t}\right),
		\end{equation}
		for some positive constants $C, c$.
		
		\smallskip
		Since
		\begin{equation*}
			\T_t\phi(x)=\int K(t,x,y)\phi(y)dy,
		\end{equation*}
		using Jensen's inequality, we obtain the following
		\begin{equation}\label{eq:Jensen_ineq_to_K}
			\begin{aligned}
				\|\nabla_x \T_t\phi\|_{L^2}^2\le& \int\left( \int |\nabla_x K(t,x,y)|\phi(y)dy\right)^2 dx\\
				=& \int\left( \int \frac{|\nabla_xK(t,x,y)|}{\int |\nabla_x K(t,x,y)|dy}\phi(y)dy\right)^2 \left(\int |\nabla_x K(t,x,y)|dy\right)^2dx\\
				\le& \int\left( \int \frac{|\nabla_xK(t,x,y)|}{\int |\nabla_x K(t,x,y)|dy}\phi^2(y)dy\right) \left(\int |\nabla_x K(t,x,y)|dy\right)^2dx\\
				\le& \int\int \Big(|\nabla_x K(t,x,y)|\phi^2(y)dy\Big)\left(\int |\nabla_x K(t,x,y)|dy\right)dx.\\
			\end{aligned}
		\end{equation}
		Now using \eqref{eq:1_diff_bound} and $u=(x-y)/\sqrt{t}$, we get the following bound
		\begin{equation}\label{eq:bound_Gaussian}
			\begin{aligned}
				\int |\nabla_x K(t,x,y)|dy \le& \int \frac{C}{t^{(d+1)/2}}\exp\left(-c\frac{|x-y|^2}{t}\right)dy\\
				\le& \frac{C}{\sqrt{t}}\int e^{-c|u|^2}du\\
				\le& \frac{C}{\sqrt{t}}\Big(\frac{\pi}{c}\Big)^{d/2}.
			\end{aligned}
		\end{equation}
		Notice that we also have the same bound for the integration with respect to $x$, i.e.
		\begin{equation}\label{eq:bound_Gaussian_x}
			\begin{aligned}
				\int |\nabla_x K(t,x,y)|dx 
				\le& \frac{C}{\sqrt{t}}\Big(\frac{\pi}{c}\Big)^{d/2}.
			\end{aligned}
		\end{equation}
		So we deduce from \eqref{eq:Jensen_ineq_to_K}, \eqref{eq:bound_Gaussian} and \eqref{eq:bound_Gaussian_x} that
		\begin{equation}\label{eq:proof:step2_case_k=0}
			\begin{aligned}
				\|\nabla_x \T_t\phi\|_{L^2}^2\le& \frac{C}{\sqrt{t}}\int\int |\nabla_x K(t,x,y)|\phi^2(y)dydx\\
				=& \frac{C}{\sqrt{t}}\int\phi^2(y)\left(\int |\nabla_x K(t,x,y)|dx\right) dy\\
				\le& \frac{C}{t}\|\phi\|_{L^2}^2.
			\end{aligned}
		\end{equation}

		\medskip
		
		\begin{itemize}
			\item \underline{Case $k=1$:}
		\end{itemize}
		
		Now to prove \eqref{eq:proof_regularity_semigroup_step2_ineq2}, we also start from equation \eqref{eq:Solution_to_parabolic_system}. By Duhamel's formula, the solution to \eqref{eq:Solution_to_parabolic_system} satisfies
		\begin{equation}\label{eq:proof_step2_Duhamel_solution}
			\psi_l(t)=\T_t\psi_l(0)+\int_0^t \T_{t-s}\divv\big((\pa_{x_l} a) \psi(s)\big)ds.
		\end{equation}
		So we have
		\begin{align}
			\|\nabla \psi_l(t)\|_{L^2}
			\le& \|\nabla \T_t\psi_l(0)\|_{L^2}+\int_0^t \|\nabla \T_{t-s}\divv\big((\pa_{x_l} a) \psi(s)\big)\|_{L^2}ds\nonumber\\
			\le& \frac{C}{\sqrt{t}}\|\psi_l(0)\|_{L^2}+C\int_0^t \frac{1}{\sqrt{t-s}}\|\divv\big((\pa_{x_l} a) \psi(s)\big)\|_{L^2}ds\nonumber\\
			\le& \frac{C}{\sqrt{t}}\|\psi_l(0)\|_{L^2}+C\int_0^t \frac{1}{\sqrt{t-s}}\|\psi(s)\|_{H^1}ds.
		\end{align}
		This implies
		\begin{align}\label{eq:proof_step2_H^1norm_ineq1}
			\sqrt{t}\|\nabla \psi_l(t)\|_{L^2}
			\le& C\|\psi_l(0)\|_{L^2}+C\int_0^t \frac{\sqrt{t}}{\sqrt{s(t-s)}}\big(\sqrt{s}\|\psi(s)\|_{H^1}\big)ds.
		\end{align}
		Now by the change of variable $s=ut$, we have
		\begin{align*}
			\int_0^t \frac{1}{\sqrt{s(t-s)}}ds=\int_0^1 \frac{1}{\sqrt{u(1-u)}}du=K<+\infty.
		\end{align*}
		Hence we can deduce from \eqref{eq:proof_step2_H^1norm_ineq1} that
		\begin{align}\label{eq:proof_step2_H^1norm_ineq2}
			\sup_{s\le t}\sqrt{s}\|\nabla \psi_l(s)\|_{L^2}
			\le& C\|\psi_l(0)\|_{L^2}+CK\sqrt{t}\sup_{s\le t}\big(\sqrt{s}\|\psi(s)\|_{H^1}\big),
		\end{align}
		and by summing up on $l=1,\ldots,d$, we have
		\begin{align}\label{eq:proof_step2_H^1norm_ineq3-1}
			\sup_{s\le t}\sqrt{s}\|\nabla\psi(s)\|_{L^2}
			\le& C\|\psi(0)\|_{L^2}+C\sqrt{t}\sup_{s\le t}\big(\sqrt{s}\|\psi(s)\|_{H^1}\big)\nonumber\\
			\le& C\|\nabla\phi\|_{L^2}+C\sqrt{t}\sup_{s\le t}\big(\sqrt{s}\|\psi(s)\|_{H^1}\big).
		\end{align}
		Now by combining \eqref{eq:proof_step2_H^1norm_ineq3-1} and \eqref{eq:proof:step2_case_k=0} we obtain
		\begin{align}\label{eq:proof_step2_H^1norm_ineq3}
			\sup_{s\le t}\sqrt{s}\|\psi(s)\|_{H^1}
			\le& C\|\phi\|_{H^1}+C\sqrt{t}\sup_{s\le t}\big(\sqrt{s}\|\psi(s)\|_{H^1}\big).
		\end{align}
		Let $t_0=1/(4C^2)$ then for $t\le t_0$, we have
		\begin{align}\label{eq:proof_step2_H^1norm_ineq4}
			\|\psi(t)\|_{H^1}\le \frac{C}{\sqrt{t}}\|\phi\|_{H^1}.
		\end{align}
		For $t> t_0$, by using the previous estimate \eqref{eq:proof_regularity_semigroup_step2_L^2_estimate_psi_Gronwall2} we also have 
		\begin{align}\label{eq:proof_step2_H^1norm_ineq5}
			\|\psi(t)\|_{H^1}\le \frac{C}{\sqrt{t_0}}\|\psi(t-t_0)\|_{L^2}\le \frac{C}{\sqrt{t_0}}\|\phi\|_{H^1}.
		\end{align}
		So gathering the two cases, we conclude that for all $t>0$,
		\begin{equation}
			\|\psi(t)\|_{H^1}\le C\left(1+\frac{1}{\sqrt{t}}\right)\|\phi\|_{H^1},
		\end{equation}
		and it implies
		\begin{align}\label{prove_step4}
			\forall\, t>0, \quad \|\nabla\T_t\phi\|_{H^1}\le& C\left(1+\frac{1}{\sqrt{t}}\right)\|\phi\|_{H^1}.
		\end{align}

		\medskip
		\paragraph{\textbf{Step 3.} }
		Now we will perform induction for \eqref{eq:semigroup_bound_first_estimate}. First we assume that \eqref{eq:semigroup_bound_first_estimate} is true for $2k$, we will prove it holds true for $2k+1$. Then, our aim is to bound $\|\nabla^{2k+1}\T_t\phi\|_{L^2}$.
		
		Using Lemma \ref{eq:Hk_estimate_for_induction}, we have the following estimate
		\begin{equation}
			\begin{aligned}
				\|\nabla^{2k+1} (\T_t\phi)\|_{L^2}=& \|\nabla^{2k} (\nabla\T_t\phi)\|_{L^2}\\
				\le& C\Big(\|A^{k}(\nabla\T_t\phi)\|_{L^2}+\|\nabla\T_t\phi\|_{H^{2k-1}}\Big)\\
				\le& C\Big(\|\nabla\T_t(A^{k}\phi)\|_{L^2}+\|[\nabla,A^k]\T_t\phi\|_{L^2}+\|\T_t\phi\|_{H^{2k}}\Big),
			\end{aligned}
		\end{equation}
		where we used the fact that $A$ commutes with $\T_t$.
		
		By the induction hypothesis, we have
		$$\|\T_t\phi\|_{H^{2k}}\le C\|\phi\|_{H^{2k}}.$$
		Since $[\nabla,A^k]$ is a differential operator of order $2k$ with smooth coefficients, we also have
		\begin{align}
			\|[\nabla,A^k]\T_t\phi\|_{L^2}\le C\|\T_t\phi\|_{H^{2k}}\le C\|\phi\|_{H^{2k}}.
		\end{align}
		Now combining with the inequality \eqref{eq:proof_regularity_semigroup_step2_L^2_estimate_psi_Gronwall2} obtained in \textbf{Step 1} case $k=1$, which allows us to treat the term $\|\nabla\T_t(A^{k}\phi)\|_{L^2}$, and then using estimate \eqref{eq:lem:Akphi_H1_estimate}, we finally get
		\begin{equation}
			\begin{aligned}
				\|\nabla^{2k+1} (\T_t\phi)\|_{L^2}
				\le& C\Big(\|A^{k}\phi\|_{H^1}+\|\phi\|_{H^{2k}}\Big)\nonumber\\
				\le& C\|\phi\|_{H^{2k+1}}.
			\end{aligned}
		\end{equation}
		
		\medskip

		Now assuming that \eqref{eq:semigroup_bound_first_estimate} is true for $2k+1$, we will prove it is true for $2k+2$. Indeed, by using Lemma \ref{lem:Hk_estimate_for_induction} and the fact that the semigroup commutes with its generator, we have
		\begin{equation}
			\begin{aligned}
				\|\nabla^{2k+2} (\T_t\phi)\|_{L^2}\le& C\Big(\|A^{k+1}(\T_t\phi)\|_{L^2}+\|\T_t\phi\|_{H^{2k+1}}\Big)\\
				\le& C\Big(\|\T_t(A^{k+1}\phi)\|_{L^2}+\|\phi\|_{H^{2k+1}}\Big)\\
				\le& C\Big(\|A^{k+1}\phi\|_{L^2}+\|\phi\|_{H^{2k+1}}\Big),
			\end{aligned}
		\end{equation}
		where we used \eqref{eq:proof_regularity_semigroup_step1_k=0} in the case $k=0$ to get the last inequality.
		
		\smallskip
		Finally, the estimate on the r.h.s. of \eqref{eq:Hk_estimate_for_induction} allows to deduce that
		\begin{equation}
			\|\nabla^{2k+2} (\T_t\phi)\|_{L^2}\le C\|\phi\|_{H^{2k+2}},
		\end{equation}
		and we complete the proof of \eqref{eq:semigroup_bound_first_estimate}.
		
		\medskip
		\paragraph{\textbf{Step 4.} } By the similar steps, we can also perform the induction steps for \eqref{eq:semigroup_bound_second_estimate} from $2k$ to $2k+1$ and from $2k+1$ to $2k+2$ respectively.
		
		First we assume that \eqref{eq:semigroup_bound_second_estimate} is true for $2k$, we prove it holds true for $2k+1$. Again we use Lemma \ref{eq:Hk_estimate_for_induction} and then apply the induction hypothesis on $2k$ to get the following
		\begin{equation}
			\begin{aligned}
				\|\nabla^{2k+1} (\nabla\T_t\phi)\|_{L^2}=& \|\nabla^{2k} (\nabla^2\T_t\phi)\|_{L^2}\\
				\le& C\Big(\|A^{k}(\nabla^2\T_t\phi)\|_{L^2}+\|\nabla^2\T_t\phi\|_{H^{2k-1}}\Big)\\
				\le& C\Big(\|\nabla^2\T_t(A^{k}\phi)\|_{L^2}+\|[\nabla^2,A^{k}]\T_t\phi\|_{L^2}+\|\nabla\T_t\phi\|_{H^{2k}}\Big).
			\end{aligned}
		\end{equation}
		Using the induction hypothesis on $2k$, we have
		\begin{equation}\label{eq:step4_1}
			\|\nabla\T_t\phi\|_{H^{2k}}\le C\left(1+\frac{1}{\sqrt{t}}\right)\|\phi\|_{H^{2k}}.
		\end{equation}
		Since $[\nabla^2,A^k]$ is a differential operator of order $2k+1$ with smooth coefficients, we also have
		\begin{equation}\label{eq:step4_2}
			\|[\nabla^2,A^k]\T_t\phi\|_{L^2}\le C\|\nabla\T_t\phi\|_{H^{2k}}\le C\left(1+\frac{1}{\sqrt{t}}\right)\|\phi\|_{H^{2k}}.
		\end{equation}
		To treat the first term, we use estimate \eqref{prove_step4} obtained in \textbf{Step 2} case $k=1$, 
		\begin{equation}\label{eq:step4_3}
			\|\nabla^2\T_t(A^{k}\phi)\|_{L^2}\le C\left(1+\frac{1}{\sqrt{t}}\right)\|A^{k}\phi\|_{H^1}
		\end{equation}
		Now combining \eqref{eq:step4_1}, \eqref{eq:step4_2} and \eqref{eq:step4_3} we obtain
		\begin{equation}
			\begin{aligned}
				\|\nabla^{2k+1} (\nabla\T_t\phi)\|_{L^2}
				\le& C\left(1+\frac{1}{\sqrt{t}}\right)\Big(\|A^{k}\phi\|_{H^1}+\|\phi\|_{H^{2k}}\Big)\nonumber\\
				\le& C\left(1+\frac{1}{\sqrt{t}}\right)\|\phi\|_{H^{2k+1}},
			\end{aligned}
		\end{equation}
		where the last inequality follows by Lemma \eqref{eq:Hk_estimate_for_induction}-$(2)$.
		
		\medskip 
		For the induction step from $2k+1$ to $2k+2$, we also start by applying Lemma \ref{eq:Hk_estimate_for_induction} and the fact that $T_t$ commutes with $A^{k+1}$ to obtain the following
		\begin{equation}
			\begin{aligned}
				\|\nabla^{2k+2} (\nabla\T_t\phi)\|_{L^2}
				\le& C\Big(\|A^{k+1}(\nabla\T_t\phi)\|_{L^2}+\|\nabla\T_t\phi\|_{H^{2k+1}}\Big)\\
				\le& C\Big(\|\nabla\T_t(A^{k+1}\phi)\|_{L^2}+\|[\nabla,A^{k+1}]\T_t\phi\|_{L^2}+\|\nabla\T_t\phi\|_{H^{2k+1}}\Big).
			\end{aligned}
		\end{equation}
		We use estimate \eqref{eq:proof:step2_case_k=0} obtained in \textbf{Step 2} case $k=0$ to bound the first term
		\begin{equation}\label{eq:step4_1'}
			\|\nabla \T_t(A^{k+1}\phi)\|_{L^2}\le \frac{C}{\sqrt{t}}\|A^{k+1}\phi\|_{L^2}\le \frac{C}{\sqrt{t}}\|\phi\|_{H^{2k+2}}.
		\end{equation}
		The induction hypothesis on $2k+1$ allows us to bound the third term
		\begin{equation}\label{eq:step4_2'}
			\|\nabla\T_t\phi\|_{H^{2k+1}}\le C\left(1+\frac{1}{\sqrt{t}}\right)\|\phi\|_{H^{2k+1}}.
		\end{equation}
		Since $[\nabla,A^{k+1}]$ is a differential operator of order $2k+2$ with smooth coefficients, we also have
		\begin{equation}\label{eq:step4_3'}
			\|[\nabla,A^{k+1}]\T_t\phi\|_{L^2}\le C\|\nabla\T_t\phi\|_{H^{2k+1}}\le C\left(1+\frac{1}{\sqrt{t}}\right)\|\phi\|_{H^{2k+1}}.
		\end{equation}
		Therefore, from \eqref{eq:step4_1'}, \eqref{eq:step4_2'} and \eqref{eq:step4_3'} we conclude that
		\begin{equation}
			\|\nabla^{2k+2} (\nabla\T_t\phi)\|_{L^2}
			\le C\left(1+\frac{1}{\sqrt{t}}\right)\|\phi\|_{H^{2k+2}}.
		\end{equation}
		The proof is completed.
	\end{proof}

	\section{The case with weights}

	Following the main steps in the proof of Proposition \ref{prop:regularity_of_semigroups}, we can obtain similar estimates for the Sobolev spaces endowed with weights. Recall that the weight $w$ is assumed to satisfy the assumptions \eqref{eq:derivative_of_w} and \eqref{eq:w(x)w(y)-condition}. The additional difficulty when working with weighted norms is that once we perform integration by parts, some extra terms are created. But thanks to the above assumption on the gradient of the weight $w$, we can bound all these extra terms by appropriate estimates.
	
	\medskip
	\begin{proof}[Proof of Threorem \ref{prop:weight:regularity_of_semigroups}]
		We follow step by step the proof of Proposition \ref{prop:regularity_of_semigroups}.
		\paragraph{\textbf{Step 1.} } We will prove \eqref{eq:weight:semigroup_bound_first_estimate} in the cases $k=0$ and $k=1$. Again we start from equation \eqref{eq:Solution_to_parabolic_eq_1}, then multiply both sides of \eqref{eq:Solution_to_parabolic_eq_1} by $\T_t\phi \, w$ and integrate in position. The only difference in comparison to \eqref{eq:proof_regularity_semigroup_step1_L^2_integration} is that: when we take integration by parts, there is an additional term created by the derivatives of the weight $w$. We have
		\begin{align}\label{eq:weight:proof_regularity_semigroup_step1_L^2_integration}
			\frac12\|\T_t\phi\|_{L^2_{w}}^2 =& \frac12\|\T_0\phi\|_{L^2_{w}}^2+\int_0^t \int_{\R^d}\T_s\phi \divv\big(a\nabla (\T_s\phi)\big)wds \nonumber\\
			=& \frac12\|\phi\|_{L^2_{w}}^2-\int_0^t \int_{\R^d}\nabla(\T_s\phi) a\nabla (\T_s\phi)wds-\int_0^t \int_{\R^d}\T_s\phi\nabla(\T_s\phi) a\nabla w ds.
		\end{align}
		Now using property \eqref{eq:derivative_of_w} of the weight $w$ we have, for any $s\in[0,t]$,
		
		\begin{align}\label{eq:weight:proof_regularity_semigroup_step1_11111}
			- \int_{\R^d}\T_s\phi\nabla(\T_s\phi) a\nabla w
			\le& C \int_{\R^d}|\T_s\phi\nabla(\T_s\phi)| w\nonumber\\
			\le& \lambda \|\nabla(\T_s\phi)\|_{L^2_{w}}^2 +\frac{C^2}{4\lambda} \|\T_s\phi\|_{L^2_{w}}^2
		\end{align}
		($\lambda$ here is the constant appearing in the uniform ellipticity condition).
		
		Now we combine \eqref{eq:weight:proof_regularity_semigroup_step1_11111} with the uniform elliptic condition \eqref{eq:proof_regularity_semigroup_step1_elliptic_cond} to deduce that
		\begin{align}\label{eq:weight:proof_regularity_semigroup_step1_k=0}
			\frac12\|\T_t\phi\|_{L^2_{w}}^2 
			\le& \frac12\|\phi\|_{L^2_{w}}^2+\frac{C^2}{4\lambda}\int_0^t \|\T_s\phi\|_{L^2_{w}}^2 ds.
		\end{align}
		Finally, using Gronwall's lemma we complete the proof of \eqref{eq:semigroup_bound_first_estimate} for the case $k=0$.
		
		\medskip
		Next, we prove \eqref{eq:weight:semigroup_bound_first_estimate} for the case $k=1$ with the same arguments as in the case without weights. Indeed, we also start by multiplying both sides of \eqref{eq:Solution_to_parabolic_system} with $\psi_l w$, then integrating in position. In the case with weight $w$, when we take integration by parts, there are also two additional terms  
		$$\int_{\R^d}\psi_l(s) \nabla w a\nabla \psi_l(s),\quad \int_{\R^d}\psi_l(s) \nabla w (\pa_{x_l} a) \psi(s).$$
		Again by using the property \eqref{eq:derivative_of_w} of weight $w$ and the fact that matrix $a$ has smooth enough coefficients, we can bound these two terms as in \eqref{eq:weight:proof_regularity_semigroup_step1_11111} and deduce the desired estimate for the case $k=1$.

		\medskip
		
		\paragraph{\textbf{Step 2.} }
		In this step, we will verify that the inequality \ref{eq:weight:semigroup_bound_second_estimate} is also true for $k=0$ and $k=1$ in the case with weights. As in \eqref{eq:Jensen_ineq_to_K}, we also have the following
		\begin{equation}
			\begin{aligned}
				\|\nabla_x \T_t\phi\|_{L^2_{w}}^2 \le& \int\left( \int \big|\nabla_x K(t,x,y)\big|\phi(y)dy\right)^2 w(x)dx\\
				=& \int\left( \int \big|\nabla_x K(t,x,y)\big|\sqrt{\frac{w(x)}{w(y)}}\phi(y)\sqrt{w(y)}dy\right)^2 dx\\
				=& \int\left( \int \big|\tilde{K}(t,x,y)\big|\tilde{\phi}(y)dy\right)^2 dx,
			\end{aligned}
		\end{equation}
		where $\tilde{K}(t,x,y)$ and $\tilde{\phi}$ are defined by 
		\begin{align}
			\tilde{K}(t,x,y)=&\nabla_x K(t,x,y)\sqrt{\frac{w(x)}{w(y)}}\\
			\tilde{\phi}(y)=&\phi(y)\sqrt{w(y)}
		\end{align}
		Next, we will prove that, for all $t\le T$,
		\begin{equation}\label{eq:weight:bound_new_kernel}
			\int \big|\tilde{K}(t,x,y)\big|dy \le \frac{C_T}{\sqrt{t}},\quad
			\int \big|\tilde{K}(t,x,y)\big|dx \le \frac{C_T}{\sqrt{t}}.
		\end{equation}
		Once the above bounds are proved, we can apply the arguments of \textbf{Step 2} case $k=0$ in the proof of Proposition \ref{prop:regularity_of_semigroups} to the kernel $\tilde{K}$ (instead of $\nabla_x K$) and function $\tilde{\phi}$ (instead of $\phi$), and obtain the desired results for the case with weights.
		
		\medskip
		
		Indeed, by using Lemma \ref{lem:differentiability_of_fundamental_solution} on the differentiability of fundamental solution and the assumption \eqref{eq:w(x)w(y)-condition} of the weight $w$, we get
		%	\begin{equation}\label{eq:proof:step2_22222}
			%		\begin{aligned}
				%			\int \big|\tilde{K}(t,x,y)\big|dy =& \int \big|\nabla_x K(t,x,y)\big|\sqrt{\frac{w(x)}{w(y)}}dy\\
				%			\le& \int \frac{C}{t^{(d+1)/2}}\exp\left(-c\frac{|x-y|^2}{t}\right)\left(\frac{1+|x|^2}{1+|y|^2} \right)^{w/2}dy\\
				%			\le& \frac{C}{t^{(d+1)/2}}\int \exp\left(-c\frac{|x-y|^2}{t}\right)\left(1+|x-y|^2 \right)^{w/2}dy.
				%		\end{aligned}
			%	\end{equation}
		%	Let $u=(x-y)/\sqrt{t}$, we can deduce that
		%	\begin{equation}\label{eq:weight:bound_Gaussian}
			%		\begin{aligned}
				%			\int \big|\tilde{K}(t,x,y)\big|dy
				%			\le& \frac{C}{\sqrt{t}}\int e^{-c|u|^2}(1+t|u|^2)^{w/2}du\\
				%			\le& \frac{C}{\sqrt{t}}\int e^{-c|u|^2}(1+t^{w/2}|u|^{w})du\\
				%			\le& \frac{C}{\sqrt{t}}\left(\sqrt{\frac{\pi}{c}}+C't^{w/2}\right)\\
				%			\le& \frac{C_T}{\sqrt{t}},
				%		\end{aligned}
			%	\end{equation}

		\begin{equation}\label{eq:proof:step2_22222}
			\begin{aligned}
				\int \big|\tilde{K}(t,x,y)\big|dy =& \int \big|\nabla_x K(t,x,y)\big|\sqrt{\frac{w(x)}{w(y)}}dy\\
				\le& \int \frac{C}{t^{(d+1)/2}}\exp\left(-c\frac{|x-y|^2}{t}\right)\exp\left(\frac{1}{2}|x-y|\right)dy.\\
			\end{aligned}
		\end{equation}
		Let $u=(x-y)/\sqrt{t}$, we can deduce that
		\begin{equation}\label{eq:weight:bound_Gaussian}
			\begin{aligned}
				\int \big|\tilde{K}(t,x,y)\big|dy
				\le& \frac{C}{\sqrt{t}}\int e^{-c|u|^2+\frac{\sqrt{t}}{2}|u|}du\\
				\le& \frac{C_T}{\sqrt{t}}\int e^{-\left(\sqrt{c}|u|-\frac{\sqrt{T}}{4\sqrt{c}}\right)^2}du\\
				\le& \frac{C_T}{\sqrt{t}},
			\end{aligned}
		\end{equation}
		where we used the fact that $\displaystyle \int e^{-\left(\sqrt{c}|u|-\frac{\sqrt{T}}{{4\sqrt{c}}}\right)^2}du=C<+\infty$ to bound in the last line.
		
		\medskip
		By the symmetry of the estimates in \eqref{eq:weight:bound_Gaussian}, we also have the same bound for $\int \big|\tilde{K}(t,x,y)\big|dx$, and hence the rest of proof of \eqref{eq:weight:semigroup_bound_second_estimate} case $k=0$ follows.
		
		\medskip
		
		Next, to prove \eqref{eq:weight:semigroup_bound_second_estimate} for case $k=1$, we start from equation \eqref{eq:proof_step2_Duhamel_solution} and notice that the difference between the case with and without weight only happens when we perform the integration by parts. In the proof of \eqref{eq:semigroup_bound_second_estimate} (the case without weight) for $k=1$, there is no appearance of integration by parts so we can follow all the previous arguments basically replacing standard $H^k$ norms by the weighted $H^k_{w}$ norms.

		\medskip
		
		\paragraph{\textbf{Step 3-4.}} As in step 3 and 4 of the proof of Proposition \eqref{prop:regularity_of_semigroups}, there is no integration by parts so we are able to follow the same lines for the case with weights. The only key point is to obtain the estimate in Lemma \ref{lem:Hk_estimate_for_induction} with weighted norms (which relies on integration by parts). We now state a weighted version of Lemma \ref{lem:Hk_estimate_for_induction} and prove it.
		
		\begin{lemma}\label{lem:weight:Hk_estimate_for_induction} 
			Let $k\ge 1$ and assume that $a_{ij}\in C^{2k-1}_b(\R^d)$.
			
			\begin{enumerate} 
				\item There exists $C_1,C_2>0$ such that
				\begin{equation}\label{eq:weight:Hk_estimate_for_induction}
					C_1\|\phi\|_{H^{2k}_{w}}\le \|A^k\phi\|_{L^2_{w}}+\|\phi\|_{H^{2k-1}_{w}}\le C_2\|\phi\|_{H^{2k}_{w}}.
				\end{equation}
				
				\item If we assume that $a_{ij}\in C^{2k}_b(\R^d)$, then there exists a positive constant $C$ such that
				\begin{equation}\label{eq:weight:lem:Akphi_H1_estimate}
					\|A^k\phi\|_{H^1_{w}}\le C\|\phi\|_{H^{2k+1}_{w}}.
				\end{equation}
			\end{enumerate}
		\end{lemma}
		
		\medskip
		\paragraph{\textit{Proof of Lemma \ref{lem:weight:Hk_estimate_for_induction}}.} The proofs of inequality on the r.h.s of \eqref{eq:Hk_estimate_for_induction} and inequality \eqref{eq:lem:Akphi_H1_estimate} do not contain integration by parts, which are only based on the expansion of $A^k\phi$ and the fact that the coefficients in this expansion are smooth enough. For this reason, we can obtain the analogous results in the weighted case without additional difficulty.
		
		\medskip
		To prove inequality on the l.h.s of \eqref{eq:weight:Hk_estimate_for_induction},
		we start by verifying that
		\begin{equation}\label{eq:weight:lowebound_A^k_techstep2}
			\begin{aligned}
				\|\phi\|_{H^2_{w}}
				\le& C\Big(\|A\phi\|_{L^2_{w}}+\|\phi\|_{H^1_{w}}\Big).
			\end{aligned}
		\end{equation}
		Indeed, thanks to \eqref{eq:derivative_of_w} and the uniform ellipticity condition, after integration by parts one has
		\begin{equation}\label{eq:weight:lowebound_A^k_techstep1}
			\begin{aligned}
				\left<\phi,A\phi\right>_{L^2_{w}}
				=&\left<\nabla \phi,a\nabla\phi\right>_{L^2_{w}}+\left<\phi\nabla w,a\nabla\phi\right>_{L^2}\\
				\ge& \lambda\|\nabla\phi\|_{L^2_{w}}^2-C\|\phi\|_{L^2_{w}}\|\nabla\phi\|_{L^2_{w}}\\
				\ge& \frac{\lambda}{2}\|\nabla\phi\|_{L^2_{w}}^2-\frac{C^2}{2\lambda}\|\phi\|_{L^2_{w}}^2.
			\end{aligned}
		\end{equation}
		Applying the above estimate to $\nabla \phi$, we get
		\begin{equation}\label{eq:weight:proof:norms_equiv:2nabla_bound}
			\begin{aligned}
				\|\nabla\nabla\phi\|_{L^2_{w}}^2
				\le& \frac{2}{\lambda}\big<\nabla \phi,A(\nabla\phi)\big>_{L^2_{w}}+\frac{C}{\lambda^2}\|\nabla\phi\|_{L^2_{w}}^2\\
				\le& \frac{2}{\lambda}\Big(\big<\nabla \phi,\nabla(A\phi)\big>_{L^2_{w}}-\big<\nabla \phi,[\nabla, A]\phi\big>_{L^2_{w}}\Big)+\frac{C}{\lambda^2}\|\nabla\phi\|_{L^2_{w}}^2\\
				\le& \frac{2}{\lambda}\Big(-\big<\Delta \phi,A\phi\big>_{L^2_{w}}-\big<\nabla\phi\nabla w,A\phi\big>_{L^2}-\big<\nabla \phi,[\nabla, A]\phi\big>_{L^2_{w}}\Big)+\frac{C}{\lambda^2}\|\nabla\phi\|_{L^2_{w}}^2.
			\end{aligned}
		\end{equation}
		Here we notice that the second term in the last line is an $L^2$-bracket (not $L^2_{w}$), which is obtained after integration by parts.
		
		Now as in \eqref{bound_for_commutator} we also have
		\begin{equation*}
			\|[\nabla, A]\phi\|_{L^2_{w}}\le C \left(\|\nabla ^2 \phi\|_{L^2_{w}}+\|\nabla\phi\|_{L^2_{w}}\right)\le C \|\phi\|_{H^2_{w}}.
		\end{equation*}
		
		By using Cauchy-Schwartz and the assumption $\nabla w\le Cw$, we deduce from \eqref{eq:weight:proof:norms_equiv:2nabla_bound} that
		\begin{equation}\label{eq:weight:bound11111111111}
			\begin{aligned}
				\|\nabla^2\phi\|_{L^2_{w}}^2
				\le& C\Big(\|\Delta\phi\|_{L^2_{w}}\|A\phi\|_{L^2_{w}}+\|\nabla\phi\|_{L^2_{w}}\|A\phi\|_{L^2_{w}}+\|\nabla \phi\|_{L^2_{w}}\|[\nabla, A]\phi\|_{L^2_{w}}+\|\nabla\phi\|_{L^2_{w}}^2\Big)\\
				\le& C\Big(\|\nabla ^2\phi\|_{L^2_{w}}\|A\phi\|_{L^2_{w}}+\|\nabla\phi\|_{L^2_{w}}\|A\phi\|_{L^2_{w}}+\|\nabla \phi\|_{L^2_{w}}\|\phi\|_{H^2_{w}}+\|\nabla\phi\|_{L^2_{w}}^2\Big).
			\end{aligned}
		\end{equation}
		
		Hence we finally obtain the following
		\begin{equation}
			\begin{aligned}
				\|\phi\|_{H^2_{w}}^2
				\le& C\Big(\|\phi\|_{H^2_{w}}\|A\phi\|_{L^2_{w}}+\|\phi\|_{H^1_{w}}\| \phi\|_{H^2_{w}}\Big).
			\end{aligned}
		\end{equation}
		In other words,
		\begin{equation}\label{weight:initial_for induction}
			\begin{aligned}
				\|\phi\|_{H^2_{w}}
				\le& C\Big(\|A\phi\|_{L^2_{w}}+\|\phi\|_{H^1_{w}}\Big).
			\end{aligned}
		\end{equation}
		
		Now, starting from \eqref{weight:initial_for induction}, we can perform the induction steps in the case with weight $w$ following the previous arguments without any additional difficulty. We just need to replace standard $H^k$ norms by the weighted $H^k_{w}$ norms.

		%\appendix
		\section{Appendix}
		\subsection{Proof of Proposition \ref{prop:weight_examples}}\label{sec:example}
		
		\begin{proof}
			$i)$
			The first condition \eqref{eq:derivative_of_w} can be easily verified. For all $\alpha \in\R_{+}$, we have
			\begin{equation}
				\begin{aligned}
					|\nabla w_{\alpha}(x)|=&2|\alpha| |x|\big(1+|x|^2\big)^{\alpha-1}\\
					\le& |\alpha| \big(1+|x|^2\big)^{\alpha}\\
					\le& |\alpha| w_{\alpha}(x).
				\end{aligned}
			\end{equation}
			
			To prove the second condition \eqref{eq:w(x)w(y)-condition}, we need the following lemma.
			
			\begin{lemma}\label{lem:ab}
				For any $x,y\in\R^d$,
				\begin{equation}\label{eq:lem:ab}
					\frac{1+|x|^2}{1+|y|^2}\le \frac{4}{3} \big(1+|x-y|^2\big).
				\end{equation}
			\end{lemma}
			
			\begin{proof}[Proof of Lemma \ref{lem:ab}.]
				By replacing $a=x-y$ and $b=y$, it is equivalent to prove that $$\frac{4}{3}\big(1+|a|^2\big)\big(1+|b|^2\big)\ge 1+|a+b|^2.$$
				Indeed, one has 
				\begin{align*}
					\frac{4}{3}\big(1+|a|^2\big)\big(1+|b|^2\big)=& 1+\frac{4}{3}\big(|a|^2+|b|^2\big)+\frac{4}{3}|a|^2|b|^2+\frac{1}{3}\\
					\ge& 1+\frac{4}{3}\big(|a|^2+|b|^2\big)+\frac{4}{3}|a||b|.
				\end{align*}
				Since $\frac{2}{3}\big(|a|^2+|b|^2\big)\ge \frac{1}{3}\big(|a|+|b|\big)^2$ and $\frac{2}{3}\big(|a|^2+|b|^2\big)+\frac{4}{3}|a||b|=\frac{2}{3}\big(|a|+|b|\big)^2$, we finally obtain
				\begin{align*}
					\frac{4}{3}\big(1+|a|^2\big)\big(1+|b|^2\big)
					\ge& 1+\big(|a|+|b|\big)^2\\
					\ge& 1+|a+b|^2.
				\end{align*}
			\end{proof}
			
			\medskip
			Now we go back to the proof of Proposition~\ref{prop:weight_examples}. Using the above lemma and combining with the inequality 
			\[2 e^x\ge1+x^2, \quad \forall\; x \in\R_{+},\]
			we conclude that 
			\begin{equation}
				\begin{aligned}
					w(x)w(y)^{-1}=&\left(\frac{1+|x|^2}{1+|y|^2}\right)^{\alpha}\\
					\le& C\left(1+|x-y|^2\right)^{|\alpha|}\\
					\le& Ce^{|\alpha||x-y|}.
				\end{aligned}
			\end{equation}

		\end{proof}
		
		$ii)$ One has
		\begin{equation}
			\begin{aligned}
				|\nabla w_{\lambda}(x)|=&\frac{|\lambda| |x|}{\sqrt{1+|x|^2}}\exp\left(\lambda\sqrt{1+|x|^2}\right)\\
				\le& |\lambda| \exp\left(\lambda\sqrt{1+|x|^2}\right)\\
				\le&|\lambda| w_{\lambda}(x),
			\end{aligned}
		\end{equation}
		so condition \eqref{eq:derivative_of_w} is satisfied.
		
		\smallskip
		Now we will verify the second condition \eqref{eq:w(x)w(y)-condition} for the weight $w_{\lambda}$. Indeed, we have
		\begin{equation}
			\begin{aligned}
				w(x)w(y)^{-1}=&\exp\left(\lambda\big(\sqrt{1+|x|^2}-\sqrt{1+|y|^2}\big)\right)\\
				=& \exp\left(\lambda\left(\frac{|x|^2-|y|^2}{\sqrt{1+|x|^2}+\sqrt{1+|y|^2}}\right)\right)\\
				\le& \exp\left(|\lambda|\left(\frac{\big||x|-|y|\big| \big(|x|+|y|\big)}{\sqrt{1+|x|^2}+\sqrt{1+|y|^2}}\right)\right)\\
				\le& \exp{\big(|\lambda||x-y|\big)}.
			\end{aligned}
		\end{equation}
		
		That completes the proof.
	\end{proof}


\begin{thebibliography}{99}
		
		
		
		
		
		\bibitem{Bertoldi-Lorenzi05} M. Bertoldi, L. Lorenzi. Estimates of the derivatives for parabolic operators with unbounded coefficients. \textit{Transactions of the American Mathematical Society}. 357. 2627-2664, 2005.
		
		\bibitem{Bismut84} J. M. Bismut. \textit{Large Deviation and Malliavin Calculus}. Progress in Mathematics 45. Birkhäuser, 1984.
		
		\bibitem{Cerrai98} S. Cerrai. Some results for second order elliptic operators having unbounded coefficients. \textit{Differential Integral Equations} 11 (4) 561 - 588, 1998.
		
		\bibitem{Cerrai01} S. Cerrai. \textit{Second Order Pde's in Finite and Infinite Dimension: A Probabilistic Approach}. Springer, 2001.
		
		\bibitem{Elworthy-Li94} K. D. Elworthy, X. M. Li. Formulae for the derivatives of heat semigroups. \textit{J. Funct. Anal.} 125, 1994.
		
		\bibitem{Friedman83} A. Friedman. \textit{Partial differential equations of parabolic type}, R.E. Krieger Pub. Co., 1983.
		
		\bibitem{Hairer15} M. Hairer, M. Hutzenthaler, A. Jentzen. Loss of regularity for Kolmogorov equations. \textit{Ann. Probab.}, 43(2):468–527, 03, 2015.
		
		\bibitem{Hauray-Pardoux-Vuong22} M. Hauray, E. Pardoux, Y. V. Vuong. Central limit theorem for a stochastic spatial epidemic model with mean field interaction, preprint, 2022.
		
		\bibitem{Hutzenthaler-Pieper21} M. Hutzenthaler, D. Pieper. Differentiability of semigroups of stochastic differential equations with Hölder-continuous diffusion coefficients, \textit{ALEA, Lat. Am. J. Probab. Math. Stat.}, 2021.
		
		\bibitem{Krylov95} N.V. Krylov, \textit{Introduction to the theory of diffusion processes}. American Mathematical Society, Providence, 1995.
		
		\bibitem{Lorenzi05} L. Lorenzi. Estimates of the derivatives for a class of parabolic degenerate operators with unbounded coefficients in $\R^N$. \textit{Annali della Scuola Normale Superiore di Pisa}, 2005.
		
		\bibitem{Lunardi97} A. Lunardi. Schauder estimates for a class of degenerate elliptic and parabolic operators with unbounded coefficients in $\R^N$. \textit{Annali della Scuola Normale Superiore di Pisa}, 1997.
		
		\bibitem{Lunardi98} A. Lunardi. Schauder theorems for linear elliptic and parabolic problems with unbounded coefficients in $\R^N$. \textit{Studia Mathematica}. 128, 1998.
		
		\bibitem{Pazy83} A. Pazy. \textit{Semigroups of linear operators and applications to partial diﬀ erential equations}, volume 44 of \textit{Applied Mathematical Sciences}. Springer-Verlag, New York, 1983.
		
		\bibitem{Prato99} G. D. Prato. Regularity results for some degenerate parabolic equation. \textit{Rivista di Matematica della Università di Parma}. Serie 6, 1999.
		
		\bibitem{Priola06} E. Priola. Formulae for the derivatives of degenerate diffusion semigroups. \textit{Journal of Evolution Equations}. 6. 577-600, 2006.
		
		\bibitem{Priola-Wang06} E. Priola, F. Y. Wang. Gradient estimates for diffusion semigroups with singular coefficients. \textit{Journal of Functional Analysis}. 236. 244-264, 2006.
		
		\bibitem{Wang97}  F. Y. Wang. On estimation of the logarithmic Sobolev constant and gradient estimates of heat semigroups. \textit{Probab. Theory Related Fields} 108, 1997.
		
	\end{thebibliography}
\end{document}